\documentclass[11pt]{article}
\usepackage{amssymb}
\newtheorem{theorem}{Theorem}
\newtheorem{lemma}{Lemma}
\newtheorem{proposition}{Proposition}
\newtheorem{corollary}{Corollary}
\title{Long Arithmetic Progressions in Critical Sets}
\author{Ernie Croot \thanks{Supported in part by an NSF grant.}}

\begin{document}

\maketitle
\begin{abstract} Given a density $0 < d \leq 1$, we show for all
sufficiently large primes $p$ that if $S \subseteq {\mathbb Z}/p{\mathbb Z}$
has the least number of three-term arithmetic progressions among all 
such sets having $\geq dp$ elements, then $S$ must contain an arithmetic
progression of length at least $\log^{1/4+o(1)} p$.  
\end{abstract}

\section{Introduction}

Given a prime $p$, we say that $S \subseteq {\mathbb Z}/p {\mathbb Z}$
is a critical set for the density $d$ if and only if 
$|S| \geq d p$ and $S$ has the
least number of three-term arithmetic progressions among all 
the subsets of ${\mathbb Z}/p{\mathbb Z}$ having
at least $d p$ elements.  In this context, a three-term arithemtic progression is
a triple of residue classes $n, n + m, n + 2m$ modulo $p$.  Note that this
includes ``trivial'' progressions, which are ones where $m \equiv 0 \pmod{p}$,
as well as ``non-trivial'' progressions, which are ones where 
$m \not \equiv 0 \pmod{p}$.  We also distinguish two different progressions,
according to how they are ordered:  The progression $n, n +m, n+2m$ is
considered different from $n+2m, n+m, n$.  

The main result
of the paper is the following theorem, which basically says that critical sets
of positive density must have long arithemtic progressions.

\begin{theorem} \label{main_theorem}  Given $0 < d \leq 1$ we have that the
following holds for all sufficiently large prime numbers $p$:  
If $S \subseteq {\mathbb Z}/p{\mathbb Z}$ is a critical set for the  
density $d$, then $S$ must contain an arithmetic progression modulo $p$
of length at least $\log^{1/4+o(1)} p$.  

Moreover, we show that for every $L \geq 1$,
$0 < d \leq 1$, and $p$ sufficiently large, there exists an arithmetic progression 
$N \subseteq {\mathbb Z}/p{\mathbb Z}$ having length at least $\log^L p$, such that
$$
|S \cap N|\ >\ |N| \left ( 1 - {1 \over \log^{1/4+o(1)}p} \right ).
$$
\end{theorem}
\bigskip

We now compare this theorem with the state-of-the-art on long progressions in
arbitrary sets of integers:  As a consequence of W. T. Gowers's deep and beautiful
proof of Szemeredi's Theorem \cite[Theorem 18.6]{gowers}, one can show that
for $0 < \delta \leq  1$, and all $x$ sufficiently large,
any set $S \subseteq \{1,2,...,x\}$ having at least $\delta x$ elements contains
an arithmetic progression of length at least $\log \log \log \log \log(x) + c(\delta)$,
for some constant $c(\delta)$.  This is a considerably shorter AP than the one
given for critical sets in our theorem above.  

There are also some results for sumsets, which give much longer AP's.  
For example, J. Bourgain \cite{bourgain} proved the interesting result 
that if $A,B \subseteq \{1,...,x\}$,
where $|A| > \delta x$, $|B| > \gamma x$, then $A+B$ has an arithmetic progression
of length at least $\exp( c (\delta \gamma \log x)^{1/3} - \log\log x)$ 
(for some $c > 0$).  I. Ruzsa \cite{ruzsa} gave an ingenious construction,
which shows that for every $0 < \epsilon  < 1/3$, and all $x$ sufficiently large, 
there exists a set $A$ having at least $b(\epsilon) x$ elements (for some 
fuction $b(\epsilon) > 0$ that depends only on $\epsilon$), such that $A+A$ has no 
arithmetic progressions longer than $\exp( \log^{2/3-\epsilon}x)$.  Then, 
B. Green \cite{green} improved Bourgain's result, and showed that a sumset
$A+B$ has an arithmetic progression of length at least
$\exp(c' (\delta \gamma \log x)^{1/2} - \log\log x)$.  We note that the length of
the progressions in these sumsets is much longer than the ones we give for 
critical sets; and so, if we could somehow prove that critical sets are sumsets 
of two large sets $A$ and $B$, then our result could possibly be improved.

There are also some impressive results on  long arithmetic progressions in
repeated sumsets $A+A+\cdots + A$ and subset sums, notably those of 
Freiman \cite{freiman}; S\' ark\H ozy \cite{sarkozy}, \cite{sarkozy2}, and \cite{sarkozy3};
Lev \cite{lev}, and \cite{lev2}; Vu and Szemeredi \cite{szemeredi} and
\cite{szemeredi2};  and J. Solymosi \cite{solymosi}. 

\bigskip

\noindent {\bf Comments:  }  The method of proof of our theorem has many 
common features with the result of B. Green \cite{green}.  In particular,
we both make use of large deviation (or concentration of measure) 
results from probability theory; and we both use techniques involving
Bohr neighborhoods.  However, the combinatorial 
aspects of our two theorems are different, which reflects the fact that 
sumsets and critical sets have different properties that must be exploited
in different ways. 

It is possible to refine the proof of our theorem, to show that critical
sets have AP's whose length depends on the density $d$; so, for example,
it might be possible to prove that critical sets 
$S \subseteq {\mathbb Z}/p {\mathbb Z}$ of density $d$ 
have a long AP for any $d > (\log\log p)^{-1}$; and, if one applies
Chang's Structure Theorem, as Green does, one can maybe get an even better
result (longer AP's holding for lower densities $d$).

\section{Proof of Theorem \ref{main_theorem}}
\bigskip

We identify $S$ with the indicator function $S(n)$, which is defined as follows:
$$
S(n)\ =\ \left \{ \begin{array}{rl} 1,\ & {\rm if\ }n \in S; \\
0,\ & {\rm otherwise.} \end{array} \right.
$$

Next, we define the discrete Fourier transform of $S(n)$ to be
$$
\hat S (a)\ =\ \sum_{0 \leq n \leq p-1} S(n) e^{2\pi i a n /p}.
$$
Then, we have that the number of 3-term arithemtic progressions in the set 
$S$ is given by
$$
\sum_{r + s \equiv 2t \pmod{p}} S(r)S(s)S(t)\ =\ 
{1 \over p} \sum_{0 \leq a \leq p-1} \hat S(a)^2 \hat S(-2a).
$$

We now write this last sum as $\Sigma_1 + \Sigma_2$, where $\Sigma_1$
is the sum over all those $a$ where 
\begin{equation} \label{large_inequality}
|\hat S(-2a)|\ >\ {p \log\log p  \over \sqrt{\log p}},
\end{equation}
and where $\Sigma_2$ is the sum over the remaining values of $a$.
From Parseval's identity we deduce the estimate
\begin{equation} \label{sigma2_estimate}
|\Sigma_2|\ \leq\ {p \log\log p \over \sqrt{\log p}} \sum_{(*)} |S(a)|^2\ 
\leq\ {d p^3 \log\log p \over \sqrt{\log p}},
\end{equation}
where the condition $(*)$ is that we sum over all $0 \leq a \leq p-1$ that do
not satisfy (\ref{large_inequality}).

We now bound the number of terms in $\Sigma_1$ from above:  
Denote this number of terms by $M$.  Then, by Parseval's identity we get that
$$
{p^2 (\log\log p)^2 \over \log p} M\ <\ 
\sum_{0 \leq a \leq p-1} |\hat S(a)|^2\ =\ d p^2,
$$
which implies
\begin{equation} \label{Sigma1_bound}
M\ <\ {d \log p \over (\log\log p)^2}.
\end{equation}

We next require the following basic lemma:

\begin{lemma} \label{dirichlet_lemma}
Suppose that $K \geq 1$, $0 \leq a_1,...,a_k \leq p-1$ and 
$$
k\ <\ {\log p \over 2K \log\log p}.
$$  
Then, for $p$ sufficiently large
there is at least one integer $1 \leq n \leq p-1$ lying in the
Bohr neighborhood defined by
\begin{equation} \label{bohr_condition}
{\rm For\ all\ }i=1,2,...,k,\ 
\left | \left | {a_i n \over p} \right | \right |\ <\ {1 \over \log^K p},
\end{equation}
where $||x||$ denote the distance from $x$ to the nearest integer.
\end{lemma}

\noindent {\bf Proof of the Lemma.  }  This is nothing more than Dirichlet's
pigeonhole argument:  We consider the $p$ vectors lying in the unit $k$-cube
$$
(a_1 y / p \pmod{1},...,a_k y / p \pmod{1}),
$$
where $y$ runs through the integers $0,1,...,p-1$.   Now, by the pigeonhole
principle, there must exist two values of $y$, say $y_1$ and $y_2$, such that
$||a_i(y_1 - y_2)/p|| < 1/\log^K p$.\ \ \ \ \ \ \ \ $\blacksquare$
\bigskip

Let $\{a_1,...,a_k\}$ be the values of $a$ satisfying 
(\ref{large_inequality}),
which are the indices of the terms in $\Sigma_1$.  
Then, we apply Lemma \ref{dirichlet_lemma} with $K = 2L$, and deduce that there is 
an integer $n_0$ satisfying (\ref{bohr_condition}).  Now, let $N$
be the arithemtic progression
$$
N\ =\ \{j n_0 \pmod{p}\ :\ 0 \leq j < \log^L p\}.
$$
We identify $N$ with its scaled indicator function
$$
N(n)\ =\ \left \{ \begin{array}{rl} {1 \over |N|},\ &{\rm if\ }n\in N;\ {\rm and} \\
0,\ & {\rm if\ }n \not \in N. 
\end{array} \right.
$$
Then, we define the Fourier transform of this scaled indicator function:
$$
\hat N(a)\ =\ \sum_{n=0}^{p-1} N(n) e^{2\pi i a n /p}.
$$

We now consider the convolution
$$
(S * N)(m)\ =\ \sum_{a+b \equiv m \pmod{p}} S(a)N(b)
\ =\ {1 \over |N|} \sum_{n \in N} S(m - n).
$$
It is obvious that 
$$
0\ \leq\ (S*N)(m) \leq 1.
$$
And, we have the following basic fact

\begin{lemma} \label{bohr_to_progression} 
Suppose that 
$$
1\ >\ \epsilon\ >\ {1 \over \log^L p}.
$$  
If $(S*N)(m) > 1 - \epsilon$, for some $0 \leq m \leq p-1$, then
$S$ contains an arithmetic progression of length at least $\epsilon^{-1}$. 
\end{lemma}

\noindent {\bf Proof of the Lemma.  }  If $(S*N)(m) > 1-\epsilon$, then we
are saying that the set $S$ contains all but $\epsilon \log^L p$ of the
residues 
\begin{equation} \label{translate_set}
m,\ m - n_0,\ m - 2n_0,\ ...,\ m - \lfloor \log^L p \rfloor n_0 \pmod{p}.
\end{equation}
Clearly, then, $S$ will contain an AP of length at least $\epsilon^{-1}$
for $\epsilon > 1/\log^L p$.\ \ \ \ \ \ \ $\blacksquare$
\bigskip

We will now show that if $S$ is a critical set, then 
\begin{equation} \label{SN_claim}
(S * N)(m)\ >\ 1 - {\log\log p \over \log^{1/4} p},
\end{equation}
for some $m$; and so, our theorem will follow from Lemma \ref{bohr_to_progression}.

To show that this is the case, suppose, for proof by contradiction, that 
(\ref{SN_claim}) fails to hold for every $0 \leq m \leq p-1$; and, let
$$
\kappa\ =\ \max \left ( 1 - {\log\log p \over \log^{1/4} p},\ 
\max_{0 \leq m \leq p-1} |(S*N)(m)| \right ).
$$
Then, define the weighting function $w(m)$ for $0 \leq m \leq p-1$ to be
$$
w(m)\ =\ \kappa^{-1} (S*N)(m).
$$
Clearly,
$$
0\ \leq\ w(m)\ \leq\ 1;
$$

Now we need the following lemma:

\begin{lemma}
Suppose that $w(m)$ is a real-valued function supported on the integers in
$[0,p-1]$, satisfying $0 \leq w(m) \leq 1$.  Then, there exists a function 
$u(m)$, also supported on the integers in $[0,p-1]$, such that 

1.  $u(m) \in \{0,1\}$ for all $m=0,1,...,p-1$;

2.  $\hat u(a)\ =\ \hat w(a) + O((\log p) \sqrt{p})$; and,

3.  $\hat u(0)\ =\ \hat w(0) + \delta,$ where $0 \leq \delta < 1$.
\end{lemma}

Before we can prove this lemma we require the following concentration
of measure result due to Hoeffding \cite{hoeffding} (also see
\cite{mcdiarmid}, Theorem 5.7):

\begin{proposition}\label{heoffding_prop}  Suppose that $v_1,...,v_r$ is a sequence of 
independent random variables where $|v_i| < 1$.  Let 
$$
\nu\ =\ E(v_1 + \cdots + v_r)\ =\ E(v_1) + \cdots +  E(v_r),
$$
and let $\Sigma = v_1 + \cdots + v_r$.  Then,
$$
P \left ( |\Sigma - \nu |\ >\ r t  \right )\ \leq\ 4\exp ( -rt^2/2).
$$
\end{proposition} 

\noindent {\bf Remark:  }  A stronger result is possible here, using Hoeffding's
theorem.  The result here is obtained as follows:  Write 
$v_i  = x_i + i y_i$, where $-1 \leq x_i, y_i \leq 1$, and then observe that
the if the ``bad event'' $|\Sigma - \nu| > rt$ occurs, then either we have the 
``bad event'' 
$|\Sigma_x - \nu_x| > r t/\sqrt{2}$ or the ``bad event'' 
$|\Sigma_y - \nu_y| > rt/\sqrt{2}$, where $\Sigma_x = x_1 + \cdots + x_r$
and $\nu_x = E(x_1 + \cdots + x_r)$, and where $\Sigma_y$ and $\nu_y$
are defined analogously.  Using Hoeffding's theorem, the probability that
either of these last two bad events occuring is at most $4 \exp( - rt^2/2)$,
as in the proposition above.
\bigskip

\noindent {\bf Proof of the Lemma.  }  
We will let $u'(m)$ be a sequence of independent 
Bernoulli random variables with distribution
$$
P(u'(m) = 1)\ =\ w(m).
$$
We note that
\begin{equation} \label{expected_u'}
E(u'(m))\ =\ w(m).
\end{equation}

Then, for each integer $a$ satisfying $0 \leq a \leq p-1$, we have that 
the Fourier transform
$$
\hat u'(a)\ =\ \sum_{j=0}^{p-1} u'(j) e^{2\pi i j a/p}
$$
can be interpreted as a sum of independent random variables as follows
$$
\hat u'(a)\ =\ v_0 + \cdots + v_{p-1},\ \ {\rm where\ \ } 
v_j = u'(j) e^{2\pi  i j a/p}.
$$

Now,
$$
E(\hat u'(a))\ =\ E(v_0) + \cdots + E(v_{p-1})\ =\ \sum_{j=0}^{p-1} E(u'(j)) e^{2\pi i j a/p}
\ =\ \hat w(a).
$$
Applying the Hoeffding proposition above, we deduce that
$$
P( |\hat u'(a) - \hat w(a)| \geq (\log p) \sqrt{p})\ <\ 4 \exp( -(\log^2 p)/2 ).
$$
Thus, the probability that
\begin{equation} \label{positive_prob_event}
{\rm For\ all\ }a=0,1,...,p-1,\ \ |\hat u'(a) - \hat w(a)|\ <\ (\log p) \sqrt{p}
\end{equation}
is at least
$$
1 - 4p\ \exp \left ( - (\log^2 p)/2 \right ),
$$
which is positive for $p \geq 11$.  

Since (\ref{positive_prob_event}) holds with positive probability, there must exist
a function $u(m)$, supported on $0,1,...,p-1$, taking the values $0$ and $1$,
and such that 
$$
{\rm For\ all\ }a=0,1,...,p-1,\ \ 
\left | \hat u (a) - \hat w(a) \right |\ <\ (\log p) \sqrt{p}.
$$
Then, by reassigning at most $O((\log p) \sqrt{p})$ of the $u(m)$'s to $0$ 
or $1$ as needed, we can get 
$$
\hat u(0)\ =\ \hat w(0) + \delta,\ 0 \leq \delta < 1,
$$
while maintaining
$$
\hat u(a)\ =\ \hat w(a) + O ((\log p) \sqrt{p})
$$
for all the other values $a=1,2,...,p-1$.  Thus, we have constructed a function
$u(m)$ which satisfies the conclusion of our lemma.\ \ \ \ $\blacksquare$
\bigskip

Now let $S'$ denote the set for which $u(m)$ is the indicator function. 
Then,  we have that 
$$
|S'|\ =\ \kappa^{-1} |S| + \delta,\ \ {\rm where\ \ } 0 \leq \delta < 1.
$$
We now estimate the number of 3AP's contained in $S'$ modulo $p$:  
This number is
\begin{eqnarray}
{1 \over p} \sum_{a = 0}^{p-1} \hat u(a)^2 \hat u(-2a)
&=&\ {1 \over p} \sum_{a=0}^{p-1} (\hat w(a) + O((\log p) \sqrt{p}) )^2 
(\hat w(-2a) + O((\log p) \sqrt{p}) )  \nonumber \\
&=&\ {1 \over p} \sum_{a=0}^{p-1} \hat w(a)^2 \hat w(-2a)\ +\ E,
\end{eqnarray}
where
$$
E\ =\ O \left ( {\log p \over \sqrt{p}} 
\sum_{a=0}^{p-1} (\log^2 p) p + (\log p)\sqrt{p} |w(a)| + 
|\hat w(a)|^2 + |\hat w(a) \hat w(-2a)| \right ).
$$
Using the Cauchy-Schwarz inequality, in combination with Parseval's
identity, one can show that
$$
E\ =\ O \left ( (\log^3 p) p \sqrt{p} \right );
$$
and so it follows that the number of 3AP's in $S'$ modulo $p$ is
\begin{eqnarray} \label{3ap_count}
&& {1 \over p} \sum_{a=0}^{p-1} \hat w(a)^2 \hat w(-2a)\ +\ O \left (  (\log^3 p) p \sqrt{p}
\right )\nonumber \\
&&\ \ \ \ \ \ \ =\ {1 \over \kappa^3 p} \sum_{a=0}^{p-1} 
\hat S(a)^2 \hat S(-2a) \hat N^2(a) \hat N(-2a)\ +\ O ((\log^3 p) p \sqrt{p} ).  
\nonumber \\
\end{eqnarray}

We now break this last sum into the two sums $\Sigma_1' + \Sigma_2'$, where
$\Sigma_1'$ is over those $0 \leq a \leq p-1$ satisfying 
(\ref{large_inequality}), and $\Sigma_2'$ is the sum for the remaining values of
$a$.   Now, for each $a$ satisfying (\ref{large_inequality}) and for each 
$n \in N$ we have from (\ref{bohr_condition}) with $K = 2L$ that
$$
\left | \left | {-2a n \over p} \right | \right |\ 
\leq\ 2 \left | \left | {an \over p} \right | \right |\ <\ {2 \over \log^L p}.
$$
for $p$ sufficiently large.  The same estimate holds for the distance from 
$an/p$ to the nearest integer.  Thus,
\begin{eqnarray}
\hat N(-2a)\ &=&\ {1 \over |N|} \sum_{n \in N} e^{2\pi i (-2a n) /p} \nonumber \\
&=&\ {1 \over |N|} \sum_{n \in N} 
\left ( 1 + O \left ( {1 \over \log^L p} \right ) \right ) \nonumber \\
&=&\ 1 + O \left ( {1 \over \log^L p} \right ); \nonumber
\end{eqnarray}
and, the same estimate holds for $\hat N(a)$.  
Thus, we conclude that 
$$
\Sigma_1'\ =\ \Sigma_1\ +\ O \left ( {p^3 \over \log^L p} \right ).
$$
We also have the estimate
$$
|\Sigma_2'|\ \leq\ {p \log\log p \over \sqrt{\log p}} \sum_{(*)} |\hat S(a)|^2
\ \leq\ {d p^3 \log\log p \over \sqrt{\log p}},
$$
where $(*)$ represents the condition that $0 \leq a \leq p-1$ such that 
$a$ does not satisfy (\ref{large_inequality}).  We note that the inequality
here follows from Parseval's identity.

Combining our estimate for $\Sigma_1'$ and $\Sigma_2'$ together with
(\ref{3ap_count}), we deduce that the number of 3AP's in $S'$ modulo $p$ is
$$
{1 \over \kappa^3 p}\left ( \Sigma_1' + \Sigma_2' \right ) + O((\log^3) p \sqrt{p})
\ =\ {1 \over \kappa^3 p} \left ( \Sigma_1 + \Sigma_2 \right ) + 
O \left ( {d p^2 \log\log p \over \kappa^3 \sqrt{\log p}} \right ).
$$
Thus, 
\begin{eqnarray} \label{3AP_wrong}
\# ( {\rm 3AP's\ in\ S' }\pmod{p})\ &=&\ {1 \over \kappa^3}\ \times\ \#( {\rm 3AP's\ 
in\ S} \pmod{p})\nonumber \\
&&\ \ \ \ \ \ \ \ \ +\ O \left ( {d p^2 \log\log p \over \kappa^3 \sqrt{\log p}} \right ).
\end{eqnarray}

We now proceed to show that this is impossible, and from our chain of reasoning
above, this would mean that (\ref{SN_claim}) holds, and therefore the theorem
would follow from Lemma \ref{bohr_to_progression}.

To show the above equation cannot hold, we require the following combinatorial
lemma, which is proved using the probabilistic method, in combination with the
second moment method:

\begin{lemma}  \label{intersect_progressions}  
Suppose $A,B \subset {\mathbb Z}/p{\mathbb Z}$ have densities $\gamma$ and
$\delta$, respectively; and, suppose that $A$ and $B$ contain 
$\alpha \gamma^3 p^2$ and $\beta \delta^3 p^2$ non-trivial 3AP's, 
respectively.  
Then, there exists a subset $C$ of ${\mathbb Z}/p{\mathbb Z}$ having 
density at least 
$$
\gamma \delta + O(p^{-1/4}),
$$
such that the number of non-trivial 3AP's lying in $C$ modulo $p$ is
at most 
$$
\alpha \beta (\gamma \delta)^3 p^2\ +\ O(p^{3/2}).
$$
\end{lemma}

\noindent {\bf Remark.  }  The same result holds if we add in trivial
AP's, since a subset $D$ of ${\mathbb Z}/p{\mathbb Z}$ can have only
$O(p)$ 3AP's, which is well within the error $O(p^{3/2})$.
\bigskip

\noindent {\bf Proof of Lemma \ref{intersect_progressions}.  }  We will find
a pair of integers $u,v$ such that $A \cap (uB + v)$ has the desired properties.
First, we show that this intersection has density very close to 
$\gamma \delta$ for almost all $0 \leq u,v \leq p-1$, by using a second moment
argument:  We suppose that $u$ and $v$ are random variables chosen independently
from $\{0,...,p-1\}$ with the uniform measure.   Then, the variance 
$V(|A \cap (uB + v)|)$ is
$$
E ( |A \cap (uB + v)|^2)\ -\ E( |A \cap (uB + v)|)^2.
$$
To compute the first expectation we express the intersection as a sum of 
indicator functions:
$$
|A \cap (uB + v)|\ =\ \sum_{b \in B} f(ub + v),
$$
where $f$ is the indicator function for the set $A$.  So, we have that
\begin{eqnarray}
E(|A \cap (uB + v)|^2)\ &=&\ \sum_{b,b' \in B} E(f(ub + v) f(ub' + v)) \nonumber \\
&=&\ {1 \over p^2} \sum_{b,b' \in B} \sum_{0 \leq u,v \leq p-1} 
f(ub+v)f(ub'+v). \nonumber
\end{eqnarray}
Now, given a pair of unequal elements $b,b' \in B$, and any two elements
$a,a' \in A$ ($a$ may equal $a'$), there is exactly one pair of numbers 
$u,v \pmod{p}$ which make $ub + v \equiv a \pmod{p}$ and $ub' + v \equiv a'
\pmod{p}$.  That is, we have that if $b'\neq b$, then there are exactly 
$|A|^2$ pairs $u,v$ which make $f(ub + v) f(ub'+v) \neq 0$ (and therefore equal
to $1$).  Thus,
$$
E(|A \cap (uB+v)|^2)\ \leq\ \gamma^2 \delta^2 p^2\ +\  |B|.
$$
The term $|B|$ comes from those pairs $b,b'$ with $b = b'$.  

To estimate $E(|A \cap (uB + v)|)$, we note that for a fixed $b \in B$
and $0 \leq u \leq p-1$, the probability that $ub + v$ lies in $A$ is 
$\gamma$.  Thus, the expected size of this intersection is 
$\gamma \delta p$.  

We now conclude that 
$$
V(|A \cap (uB + v)|)\ \leq\ |B|\ =\ \delta p;
$$
and so, by an application of Chebychev's inequality we conclude that
$$
P( |A \cap (uB + V)| < (1-\epsilon ) \gamma \delta p)\ \leq\ 
{1 \over \epsilon^2 \gamma^2 \delta p}.
$$ 

Next, we compute the expected number of 3AP's in the intersection 
$A \cap (uB  + v)$:  Let $Q = Q(u,v)$ be the number of non-trivial 3AP's
lying in $A \cap (uB + v)$.  Now, 
suppoe that $x_1,x_2,x_3$ is a non-trivial 3AP in $A$, so that
$x_2 \equiv x_1 + d, x_3 \equiv x_1 + 2d \pmod{p}$, for some 
$d \not \equiv 0 \pmod{p}$;
and, suppose that $y_1,y_2,y_3$ is a non-trivial 3AP in $B$.  Then, there is exactly 
one pair $0 \leq u,v \leq p-1$ such that 
$$
{\rm For\ }i=1,2,3,\ \ ux_i + v \equiv y_i \pmod{p}.
$$
Thus, for $u,v$ chosen at random from $0,...,p-1$ with uniform probability, the
probability that a particular non-trivial 3AP lies in $uB + v$ is 
$\beta \delta^3$; and so, the expected size of $Q$ 
is $\alpha \beta (\gamma \delta)^3 p^2$.
So, there can be at most $p^2 - p^{3/2}$ of the choices for 
$u$ and $v$ such that the intersection has more than 
$\alpha \beta (\gamma \delta)^3 (p^2 + 2p^{3/2})$ 
3AP's; else, if all but $p^{3/2}$ of the choices give more than this many
3AP's in this interesection, then we would have that $Q$ exceeds
$$
{(p^2 - p^{3/2}) (p^2 + 2p^{3/2}) \over p^2} \alpha \beta (\gamma \delta)^3,
$$
which we know is not the case.
Thus, the probability that $Q < \alpha \beta (\gamma \delta)^3
(p^2 + 2p^{3/2})$ is $> p^{-1/2}$.  
So, for $\epsilon = p^{-1/4}\gamma^{-1}\delta^{-1/2}$,
we get that the events
$$
|A \cap (uB + v)|\ \geq\ (1-\epsilon) \gamma \delta p\ \ {\rm and\ \ } 
Q\ <\ \alpha \beta (\gamma \delta )^3 (p^2 + 2p^{3/2})
$$
occur with positive probability.  So, there is a choice for $u$ and $v$ 
so that both these events occur, which proves the lemma.\ \ \ \ $\blacksquare$

We require one more lemma before we can prove that (\ref{3AP_wrong}) is
impossible:

\begin{lemma}  \label{small_progressions} Given $0 < \theta < 1$, 
there exists a subset $U \subset {\mathbb Z}/p{\mathbb Z}$ having density
$1-\theta + O(1/p)$ such that the number of 3AP's lying in $U$, both trivial and 
non-trivial, is at most 
$$
p^2 ( 1 - 3\theta + 2.5\theta^2).
$$
For $0 < \theta < 1/3$ this quantity is at most 
$$
p^2 (1-\theta)^3 (1 - \theta^2/2).
$$
\end{lemma}

\noindent {\bf Proof.  }  First, we claim that the sum of the number of 3AP's
(trivial and non-trivial) lying in $U$ and lying in 
$\overline{U} = ({\mathbb Z}/p{\mathbb Z}) \setminus U$ is 
\begin{equation} \label{3AP_count}
p^2 (1 - 3\theta + 3\theta^2).
\end{equation}
This follows by inclusion-exclusion:  The number of 3AP's lying in $U$ is 
$x_1 - x_2 + x_3 - x_4$, where $x_1$ is the total number of 3AP's among
the residue classes modulo $p^2$; $x_2$ is the sum of the 
number of these 3AP's
$x,x+d,x+2d$ such that $x \in \overline{U}$, summed with the number where
$x+d \in \overline{U}$, and summed with the number where $x+2d \in \overline{U}$; 
$x_3$ is the sum of the number of such progressions with $x,x+d \in \overline{U}$, 
then $x,x+2d \in \overline{U}$, and finally summed with the number where 
$x+d,x+2d \in \overline{U}$;
finally, $x_4$ is the number of progressions in $\overline{U}$.  It is easy to see that
$x_1 = p^2$, $x_2 = 3rp^2$, and $x_3 = 3r^2 p^2$.  So,  sum of the number of 3AP's
in $U$ and $\overline{U}$ equals the expression in (\ref{3AP_count}).

Now consider the set $\overline{U}$, having density $\theta + O(1/p)$, 
given as follows:
$$
\overline{U}\  :=\ [0,\theta p/2]\ \cup\ [p/2, p/2 + \theta p/2],
$$
where here we take the integers in these two intervals (since $U$ and 
$\overline{U}$ are sets of residue classes modulo $p$).  Call the first interval
$I_0$, and the second $I_1$.  If $x,y \in I_i$, then $z \equiv 2^{-1} (x+y) \pmod{p}$
lies in $I_i$ if $x,y$ have the same parity, and lies in $I_{1-i}$ if they are
of different parity.  This gives that the number of 3AP's $x,x+d,x+2d$ is at 
least the number of ordered pairs $x,y$ where both $x,y \in I_0$ or both are in $I_1$.
So, the number of 3AP's in $\overline{U}$ is at least $\theta^2p^2/2$, and it
follows that the number of 3AP's in $U$ is at most 
$$
p^2 ( 1 - 3\theta + 3\theta^2 - \theta)^2/2)\ =\ p^2 (1 - 3\theta + 2.5\theta^2),
$$
which proves the lemma.\ \ \ \  $\blacksquare$
\bigskip

Now we let $\theta = 1 - \kappa$, and let $U$ be the set given by this lemma.  
Then, we apply Lemma \ref{intersect_progressions} with 
$A = U$, and $B = S'$, and we deduce that there is a set $C$ with
$$
|C|\ =\ |S|  + O(p^{3/4}),
$$
such that $C$ contains at most
\begin{eqnarray}
&& \kappa^3 \left ( 1 - {(1-\kappa)^2 \over 2} \right ) \left ( {\#( {\rm 3AP's\ in\ }S) \over 
\kappa^3} + O \left ( {\log\log p \over \sqrt{\log p}} \right ) \right ) \nonumber \\
&&\ \ \ \ \ \ \ \ \ \ \ =\ \#({\rm 3AP's\ in\ }S ) \left ( 1 - {(\log\log p)^2 \over 2 \sqrt{\log p}} 
+ O \left ( {\log\log p \over \sqrt{\log p}}  \right ) \right ).
\nonumber 
\end{eqnarray}

To show that this is impossible for sufficiently large $p$, we let $C'$ be any set
gotten from $C$ by adding or removing at most $O(p^{3/4})$ elements such that 
$$
|C'|\ =\ |S|.
$$
Then, in the worst case, each element we add to $C$ (to produce $C'$)  adds
at most $p$ new 3AP's.  Thus,
\begin{eqnarray} \label{C_inequality}
\#(3AP's\ {\rm in\ }C')\ &=&\ \#(3AP's\ {\rm in\ }C)\ -\ O( p^{1.75} )\nonumber \\
&<&\ \#( {\rm 3AP's\ in\ } S ) 
\left ( 1 - { (\log\log p)^2 \over 
2 \sqrt{\log p}} + 
O \left ( {\log\log p \over \sqrt{\log p}} \right ) \right ). \nonumber \\
\end{eqnarray}
To get this inequality we have used a corollary of the following theorem of Varnavides 
\cite{varnavides}, which allows us to absorb the error term $O(p^{1.75})$ into
the error $O( (\log\log p)/\sqrt{\log p})$:

\begin{theorem} Given $0 < \alpha  \leq 1$, there exists $0 < c \leq 1$ such that
for any set $T \subseteq \{1,2,...,x\}$ having $|T| \geq \alpha x$, 
$$
\#(a,b,c \in T\ :\ a+b  =2c)\ >\ cx^2.
$$
\end{theorem}

This corollary is: 

\begin{corollary} There exists $0 < c \leq 1$, depending only on $d$ (the
lower bound for the density of $S$), such that 
$$
\#( {\rm 3AP's\ in\ }S)\ >\ c p^2.
$$
\end{corollary}

The proof of this corollay is immediate, since if we think of $S$ as a set of
integers, say $S \subseteq \{0,1,...,p-1\}$ (instead of as a set of residue classes
modulo $p$), then every solution to $a+b = 2c$, $a,b,c \in S$ in the integers
gives a solution $a+b \equiv 2c \pmod{p}$.  So, the number of 3AP's in 
$S$ modulo $p$ is at least the number of 3AP's in $S$, when we think of it as a
subset of the integers.
\bigskip

Now, (\ref{C_inequality}) contradicts the fact that $S$ is a critical set:  
Here we have constructed a set
$C'$ having the same cardinality as the set $S$, but where $C'$ has fewer
3AP's than $S$.  Thus, we must conclude that
$$
c(n)\ >\ 1 - {\log\log p \over \log^{1/4} p}  
$$
for some $0 \leq n \leq p-1$, and the theorem is proved.\ \ \ \ $\blacksquare$

\end{document}